# Harmonic Center of an n-Simplex: Some Properties

V. S. Patwardhan[1]

02 May 2022

**Abstract**

A simplex in n dimensions is defined by the usual linear inequality constraints specified as **Ax** ≤ **b**, where **A** is a [(n+1) x n] matrix. Here we consider simplexes which are bounded sets. The harmonic center has been defined earlier for polytopes in general. A relationship between the rows of **A** is derived here, which must be satisfied for the simplex to be bounded. Using this condition, an interesting relationship is derived between the (n+1) residuals at the harmonic center of the simplex. Finally, a linear expression is developed which is invariant everywhere in the simplex.

## 1. Introduction

The calculation of the center of a polytope has been of interest for a long time. In this paper the term polytope is used to indicate a bounded convex set in n dimensions, defined as the intersection of m half spaces, and is assumed to be non-empty. Many different definitions have been used for the center of a polytope, which include the center of mass, the centroid, the center of the largest inscribed sphere (or an ellipsoid), the center of the smallest sphere (or ellipsoid) which includes the polytope, the analytical center, the weighted projection center, orthogonal projections on to polytope faces, the BI center, and the harmonic center [Moretti, 2003; Patwardhan, 2019, 2020]. All these definitions lead to different points as centers, whose calculation involves different degrees of computational effort. Some of these, and a few others such as the Fermat-Torricelli point, orthocentre and the Monge point, have been investigated in detail for simplexes [Edmonds et.al., 2004].

In this paper, we consider a bounded simplex in n dimensions, defined by the usual set of constraints **Ax** ≤ **b**, where **A** is a [(n+1) x n] matrix, with rank = n. It is assumed here that all the (n+1) constraints are inequality constraints. The harmonic center has been defined earlier for a polytope. Here we show that the **A** matrix must satisfy a specific condition for the simplex to be a bounded set. We then derive an interesting relationship between the residuals for all the constraints, calculated at the harmonic center of the simplex. Finally, we present a linear expression which has a constant value everywhere in the simplex.

## 2. A necessary condition for the boundedness of a simplex

We consider here, simplexes which are bounded. Any general set of (n+1) half spaces (i.e. constraints) may not give a bounded set as it may contain a ray originating from a feasible point, which can be extended indefinitely without violating any constraint. Let us consider this aspect further.

---

[1] Independent researcher. Formerly, Scientist G, National Chemical Laboratory, Pune 411008, India.
  Email : vspatw@gmail.com , URL : https://www.vspatwardhan.com

For a simplex in n dimensions defined by (n+1) inequality constraints, the $i^{th}$ constraint can be written as

$$\sum_{j=1}^{n} A_{ij} x_j \leq b_i \quad for\ i = 1\ to\ (n+1) \tag{1}$$

It is assumed here that the rows of **A** are normalized. In other words, $\sum_{j=1}^{j=n} A_{ij}^2 = 1\ for\ all\ i$. Let P be a strictly feasible, interior point. The residuals at point P can be denoted as $S_i$, given by

$$S_i = b_i - \sum_{j=1}^{n} A_{ij} P_j \quad for\ i = 1\ to\ (n+1) \tag{2}$$

Since P is an interior point, $S_i > 0$ for all i. Since the rank of A is n, we can express the last row of A as a linear combination of rows from 1 to n, in the following form:

$$\sum_{i=1}^{n} \gamma_i A_{ij} + A_{n+1,j} = 0 \quad for\ all\ j = 1\ to\ n \tag{3}$$

where $\gamma_i$ is a row multiplier. The residual for the last constraint can now be written as

$$S_{n+1} = b_{n+1} + \sum_{j=1}^{n} \left( \sum_{i=1}^{n} \gamma_i A_{ij} \right) P_j \tag{4}$$

$$= b_{n+1} + \sum_{i=1}^{n} \gamma_i \left( \sum_{j=1}^{n} A_{ij} P_j \right) \tag{5}$$

$$= b_{n+1} + \sum_{i=1}^{n} \gamma_i (b_i - S_i) \tag{6}$$

$$= b_{n+1} + b_0 - \sum_{i=1}^{n} \gamma_i S_i \tag{7}$$

where $b_0$ is given by

$$b_0 = \sum_{i=1}^{n} \gamma_i b_i \tag{8}$$

Since $S_{n+1} > 0$ for the interior point P, we get

$$b_{n+1} \geq -b_0 + \sum_{i=1}^{n} \gamma_i S_i \tag{9}$$

Let L be the straight line drawn from P in the direction given by an arbitrary unit vector **v**. Let Q be a point where this line intersects the $i^{th}$ constraint. The distance between P and Q is given by

$$d_i = \frac{S_i}{\sum_{j=1}^{n} A_{ij} v_j} \qquad for\ all\ i = 1\ to\ (n+1) \tag{10}$$

Writing this equation for i = (n+1), substituting for $A_{n+1,j}$ using equation (3) and rearranging terms algebraically, we get

$$d_{n+1} = -\frac{S_{n+1}}{\sum_{i=1}^{n} \gamma_i \beta_i} \tag{11}$$

where

$$\beta_i = \sum_{j=1}^{n} A_{ij} v_j = \frac{S_i}{d_i} \qquad for\ all\ i = 1\ to\ n \tag{12}$$

It may be noted that equation (12) can be written as $(\mathbf{A}^N)\mathbf{v} = \boldsymbol{\beta}$, where $\mathbf{A}^N$ is the square matrix corresponding to rows 1 to n of the matrix **A**. Since $\mathbf{A}^N$ is a non-singular matrix, we can calculate **v** for any given **β**.

Let us now consider the sign of $\gamma_k$. Let us choose **β** such that

$$\beta_i = -h \qquad for\ i = k \tag{13}$$

$$= 0 \qquad otherwise \tag{14}$$

where h is an arbitrary positive number. The corresponding direction **v** can be calculated as stated above, if required. Since $\beta_i = S_i / d_i$, the line in this **v** direction does not intersect any of the first n constraints. (Specifically, it is parallel to all n constraints, except constraint k, which it intersects only if extended backwards.) For the **β** selected, equation (10) simplifies to

$$d_{n+1} = -\frac{S_{n+1}}{\gamma_k \beta_k} \tag{15}$$

Since the simplex is assumed to be bounded, this line must intersect the $(n+1)^{th}$ constraint, i.e. $d_{n+1}$ must be positive. Therefore $\gamma_k > 0$. This is valid for all k = 1 to n. Therefore,

$$\gamma_k > 0 \qquad for\ all\ k = 1\ to\ n \qquad (16)$$

This proves that for a simplex to be bounded, all the γ coefficients appearing in equation (3) must be strictly positive.

### 3. A relationship between residuals at the harmonic center of a simplex

The harmonic center H for a polytope has been defined earlier (Patwardhan, 2019), and was essentially inspired by geometrical considerations. It can be briefly described as follows. If L is a straight line drawn in any direction from any feasible point P, then in general, it intersects every constraint at one point, giving m intersections, for a polytope defined by m constraints. It was shown that there exists a unique feasible point Q somewhere along this line, such that $\sum_{i=1}^{i=m}(1/d_i) = 0$, where $d_i$ is the algebraic distance between Q and the intersection with the $i^{th}$ constraint, measured along the line L. The point Q is defined as the harmonic point along the line L. The harmonic center H of the polytope is defined as that point which is the harmonic point for all n lines drawn through it, each parallel to one of the coordinate axes. It was also shown that H is in fact, the harmonic point for any line drawn through H in any direction. It was later pointed out that the harmonic center H coincides with the analytic center [Izmestiev I., 2020]. The analytic center is well-known in the field of optimization and minimises the logarithmic potential function over the feasible space of a polytope. Details of the analytic center can be found in numerous papers in literature [for example, Wang et.al., 2019]. However, we refer to it as the harmonic center here, since we are interested in bringing out some properties of this center which are geometrical in nature. Here we consider only simplexes, as mentioned earlier.

Let P be a strictly feasible, interior point inside the simplex. A line drawn through P, parallel to the $k^{th}$ axis, intersects the $i^{th}$ constraint after travelling a distance $d_{ik}$ given by

$$d_{ik} = \frac{S_i}{A_{ik}} \qquad for\ all\ k = 1\ to\ n \qquad (17)$$

Since the rows of A are assumed to be normalized, $d_{ik}$ is the actual distance. At the harmonic center H of the simplex, we can write

$$\sum_{i=1}^{n+1} \frac{1}{d_{ik}} = 0 \qquad for\ all\ k = 1\ to\ n \qquad (18)$$

In other words,

$$\sum_{i=1}^{n} \frac{A_{ik}}{S_i^H} + \frac{A_{n+1,k}}{S_k^H} = 0 \qquad for\ all\ k = 1\ to\ n \qquad (19)$$

where the superscript H indicates the residuals at the harmonic center H. Using Equation (3), this can be rewritten as

$$\sum_{i=1}^{n} A_{ik} \left[ \frac{1}{S_i^H} - \frac{\gamma_i}{S_{n+1}^H} \right] = 0 \qquad \text{for all } k = 1 \text{ to } n \tag{20}$$

Equation (20) is really a set of n equations, and can be written in a matrix form as

$$(A^N)^T f = 0 \tag{21}$$

where the superscript T indicates a transpose, and **f** is an n-vector with

$$f_i = \frac{1}{S_i^H} - \frac{\gamma_i}{S_{n+1}^H} \qquad \text{for all } i = 1 \text{ to } n \tag{22}$$

Since **A^N** is a non-singular matrix, Equation (21) is satisfied only by **f** = 0. This gives

$$\frac{1}{S_i^H} - \frac{\gamma_i}{S_{n+1}^H} = 0 \qquad \text{for all } i = 1 \text{ to } n \tag{23}$$

In other words,

$$\gamma_i S_i^H = S_{n+1}^H \qquad \text{for all } i = 1 \text{ to } n \tag{24}$$

This is an interesting relationship between residuals at the harmonic center H and can be used as a termination criterion for any iterative approach for determining H for a simplex.

### 4. An invariant linear expression for a simplex

It is trivial to show that the sum of the lengths of three orthogonals from any point inside an equilateral triangle is constant. We can use equation (7) to make a similar statement for the residuals for a simplex. For any feasible point P in the interior of a simplex, equation (7) gives

$$S_{n+1} + \sum_{i=1}^{n} \gamma_i S_i = b_{n+1} + b_0 = constant \tag{25}$$

This is a linear expression involving residuals for a simplex, which has a constant value at every point P inside the simplex.

### 5. The value of residuals at the harmonic center of a simplex

Equation (25) is valid at any feasible point. Writing this equation for the harmonic center H, and using equation (24), we can write

$$S_{n+1}^H = \frac{(b_{n+1} + b_0)}{(n+1)} \tag{26}$$

and

$$S_i^H = \frac{(b_{n+1} + b_0)}{\gamma_i(n+1)} \qquad \text{for all } i = 1 \text{ to } n \qquad (27)$$

where $b_0$ is given by equation (8).

## 6. Simplexes with equal $\gamma_i$ values

Let us consider a simplex where the last constraint has been obtained by adding up all other constraints, and changing the right-hand side of this constraint (otherwise, all constraints would pass through a common point, and the simplex would get reduced to one point), and normalizing it. Such a simplex would be characterised by equal values of $\gamma_i$ (= $\gamma$) for i = 1 to n. $\gamma$ can be calculated using equation (3), which gives

$$A_{n+1,j} = -\gamma \sum_{i=1}^{n} A_{ij} \qquad \text{for all } j = 1 \text{ to } n \qquad (28)$$

Since the last constraint is normalised, we have $\sum_{j=1}^{j=n}(A_{n+1,j}^2) = 1$. Using this, we get

$$\gamma = 1 \Big/ \left[ \sum_{j=1}^{n} \left( \sum_{i=1}^{n} A_{ij} \right)^2 \right]^{1/2} \qquad (29)$$

Equation (27) implies that all $S_i$ values for i = 1 to n, are equal at the harmonic center. In other words, the harmonic center is equidistant from constraints i = 1 to n. Its distance from the (n+1)$^{th}$ constraint is of course different and is equal to $S_{n+1}$, given by equation (26).

## 7. Conclusions

The harmonic center has been defined earlier for a closed convex polytope. Here, some interesting properties of the harmonic center for simplexes in n-dimensions, have been considered. A necessary condition for a simplex to be bounded, is derived. The residuals at the harmonic center of a simplex obey the relationship given by equation (24). A linear expression involving residuals has been derived which is invariant at all points inside the simplex. The residuals at the harmonic center of a simplex are given by equations (26) and (27).